\documentclass[a4paper,oneside,11pt]{article}
\usepackage[english]{babel}
\usepackage[latin1]{inputenc}
\usepackage{indentfirst}
\usepackage{amssymb}
\usepackage{amsmath}
\usepackage{amsthm}

\Large
\theoremstyle{plain}
 \newtheoremstyle{miestilo}{12pt}{\topsep}{\itshape}{}{\bf}{}{ }{}
  \swapnumbers
 \theoremstyle{miestilo}
\newtheorem{teorema}[subsection]{Theorem.}
\newtheorem{proposicion}[subsection]{Proposition.}
\newtheorem{lema}[subsection]{Lemma.}
\newtheorem{corolario}[subsection]{Corollary.}

 \newtheoremstyle{misnotas}{12pt}{12pt}{}{}{\bf}{}{ }{\remark}
  \theoremstyle{misnotas}
 \newtheorem{remark}[subsection]{\ {\bf Remark.}}
  \newtheorem{apartado}[subsection]{\ {\ }}
 \allowdisplaybreaks[2]

\begin{document}

\flushbottom
%
%

{\LARGE  {\centerline{PI theory for Associative Pairs} }} %

%
%
   %
 \medskip
\medskip
\centerline{ Fernando Montaner \footnote{Partially supported   by grant MTM2017-83506-C2-1-P (AEI/FEDER, UE), and by
Diputaci\'on General de Arag\'on (Grupo de Investigaci\'{o}n
\'{A}lgebra y Geometr\'{\i}a).}} \centerline{{\sl Departamento de Matem\'{a}ticas,
Universidad de Zaragoza}} \centerline{{\sl 50009 Zaragoza,
Spain}} \centerline{E-mail: fmontane@unizar.es} \centerline{and}
\centerline{ Irene Paniello \footnote{Partially supported   by grant MTM2017-83506-C2-1-P (AEI/FEDER, UE).}} \centerline{{\sl
Departamento de Estad\'{\i}stica, Inform\'{a}tica y  Matem\'{a}ticas, Universidad
P\'{u}blica de Navarra}}
 \centerline{{\sl 31006 Pamplona, Spain}}
\centerline{E-mail: irene.paniello@unavarra.es}

%
%
\begin{abstract}
We extend the classical associative PI-theory to Associative Pairs, and in doing so, we introduce related notions already present for algebras (and Jordan systems) as the ones of PI-element and PI-ideal, extended centroid and central closure.

\end{abstract}

\section*{Introduction}

Associative Pairs are the natural generalization of associative algebras in the context of Jordan pairs, which in turn, are related to Jordan and associative triple systems. As associative algebras do in the theory of Jordan algebras, they play a central role in the theory of Jordan Pairs. As a consequence, most of the questions on associative algebras that arise from Jordan theory, are also of great interest in the case of Associative Pairs (this is what could be named ``generalized Herstein Theory" after the line of research inaugurated by Herstein on the Jordan (and Lie) structures of associative algebras). In particular, that is the case with the theory of polynomial identities (PI-theory), which plays a central role in modern Jordan theory, after the sweeping work of Zelmanov on the subject.

In the present paper, we address that area for the case of Associative Pairs. Our objective is to extend most of the classical associative PI-theory. Namely, we deal with the classical theorems of that theory that are related to the structure theory, so to the notions of simplicity, primitivity or primeness, and their versions under the presence of an involution (again, a central ingredient for the associative systems appearing in Jordan -and Lie-theory) due to Kaplansky, Posner, Amitsur and Martindale, among many other authors.

There are several features which are peculiar to Pairs in contrast to algebras. Apart from the structure theoretic particularities (often advantages, as in the study of the socle and the use of idempotents), from the viewpoint of PI-theory a feature which is always present is the possibility of moving to the theory of algebras by means of the so called standard imbeddings. Although this is often very useful (and, even necessary), we have preferred whenever possible, to adapt the algebra notions to Pairs (this is due to the fact that for the applications to Jordan theory we do not dispose of a construction as the standard embedding, so using that for associative pairs arising in the Jordan context  involves a detour which may obfuscate the arguments involved). On the other hand, some ideas which are natural in the Pair context, may shed light on developments of the associative theory. In the case of the theory of polynomial  identities, this is quite notably the situation in the study of generalized polynomial identities as presented in \cite{Rowen proceedings}.

Our study of polynomial identities on Associative Pairs has a peculiarity that deserves to be stressed. We do not really work with polynomial identities as such, although this could be done by resorting to the enveloping algebras. Instead, we adopt a more pair-theoretical approach as it is the use of homotope polynomial identities, and of PI-elements, in the line of what was done \cite{pi-i, pi-ii, hpi} for Jordan systems.

This means, among other this, that we need to introduce the notion of extended centroid of an associative pair, and the corresponding scalar extension, its central closure. This raises many open problems that we do not address in the present paper, but whose solution would be doubtless quite interesting.

This paper is organized as follows.
After this introductory section, in the first section we settle the basic notation,  and recall some fundamental results on associative pairs mainly dealing with  the relationship between associative pairs and their associative standard imbeddings. In the second section, we extend the construction of the  extended centroid for semiprime associative algebras to semiprime associative pairs, following the approach of \cite{amitsur sprime rings,BM-centroid-sprime,martindale prime rings}.  Accordingly, elements of the
 extended centroid of a semiprime associative pair are defined as equivalence classes of partially defined pair homomorphims (that is pair homomorphims defined over essential ideals, and not on the whole pair) commuting with all left, right and middle multiplication operators defined by pair elements.  As proved in that section, and in a way similar to the case of associative algebras, the extended centroid of a semiprime associative pair is a commutative, unital, von Neumann regular ring, which in fact, is isomorphic to the extended centroid of the standard imbedding of the associative pair, so to the extended centroid of an associative algebra.

 In the third section we consider the central closure of semiprime associative pairs, that is the  natural scalar extension  associated to the extended centroid, which is a tight scalar extension of the associative pair, and whose standard imbedding turns out to be isomorphic to the central closure of the  standard imbedding
 of the original associative pair.

 In section  fourth we examine semiprime associative pairs endowed with polarized involutions. Such involutions have straightforward extentions  to both the extended centroid and the central closure, so allowing the study of the $\ast$-extended centroid of a semiprime associative pair, that is the set of symmetric elements of the extended centroid under the extended involution, and of a new scalar extension, the $\ast$-central closure, that is the scalar extension linked to the  $\ast$-extended centroid. Again these two constructions behave well with standard imbeddings, and it is not difficult to prove results analogous to those contained in the previous sections relating the  $\ast$-extended centroid  and the $\ast$-central closure of a semiprime associative pair with involution to those of their standard imbeddings.

 Finally,  in the fifth section, we deal with what is the central objective of the paper, namely the study of  prime and primitive associative pairs having nonzero local algebras which satisfy  polynomial identities.  We introduce the notion of strongly primitive associative pair following \cite{Rowen proceedings,Rowen PI ring}, to be an
 associative pair  with nonzero socle which is a  dense subpair  of pairs of homomorphims between two right vector spaces   over a division PI-ring and show that the strong  primitivity  of an associative pair is equivalent to that of its standard imbedding. Then analogous results to  Amitsur, Kaplansky, Martindale and Posner Theorems are given for associative pairs, based on the existence of either local PI-algebras or on the fact that the associative pair satisfies some homotope polynomial identity.

\section{Preliminaries}

\begin{apartado} We will work with associative systems (algebras and pairs) over a unital commutative ring of scalars $\Phi$ that will  be fixed throughout. We refer to \cite{loos-jp} and \cite{meyberg} for notation, terminology and basic results. In this section, we recall some of those basic notation and results.
\end{apartado}

\begin{apartado}
We denote operations $A^\sigma\times A^{-\sigma}\times A^\sigma \rightarrow A^\sigma$ of associative pairs $A = (A^+, A^-)$ over $\Phi$ by juxtaposition: $(x^\sigma, y^{-\sigma}, z^\sigma)\mapsto  x^\sigma y^{-\sigma}  z^\sigma$. We will also make use of the operator notation: $xyz=L(x,y)z=R(y,z)x=M(x,z)y,$
where $L$, $R$ and $M$ are the left, right and middle multiplication operators respectively.
\end{apartado}

  \begin{apartado}\label{subsect standard imbedding}
 For any unital associative algebra ${\cal E}$ with an idempotent $e$, its  associated Peirce decomposition ${\cal E}={\cal E}_{11}\oplus{\cal E}_{12}\oplus{\cal E}_{21}\oplus{\cal E}_{22}$   gives rise to the associative pair $A=({\cal E}_{12},{\cal E}_{21})$ with operations inherited from the multiplication in ${\cal E}$.

 Reciprocally, associative pairs are abstract off-diagonal Peirce spaces of associative algebras  \cite[p.~92, p.~101]{loos-jp}: given an associative pair  $A=(A^+,A^-)$ we can construct a unital associative algebra  ${\cal E}$ with a Peirce decomposition  ${\cal E}={\cal E}_{11}\oplus{\cal E}_{12}\oplus{\cal E}_{21}\oplus{\cal E}_{22}$, where $A=({\cal E}_{12},{\cal E}_{21})$. The $\Phi$-module ${\cal E}_{ii}$ for $i=1,2$ is the subalgebra of $End_\Phi(A^\sigma)\times End_\Phi(A^{-\sigma})^{op}$, where $\sigma =+$ if $i =1$, and $\sigma =-$ if $i=2$, generated by the idempotent $e_i=(Id_{A^\sigma}, Id_{A^{-\sigma}})$ (hence $e_1+e_2 = 1$), and all  elements $x^\sigma y^{-\sigma}=(L(x^\sigma,y^{-\sigma}), R(y^{-\sigma},x^\sigma))$.

 It is clear then  that $A^+={\cal E}_{12}$ is an ${\cal E}_{11}-{\cal E}_{22}$ bimodule, and $A^-={\cal E}_{21}$ is an ${\cal E}_{22}-{\cal E}_{11}$ bimodule  with the obvious actions (so that, in fact, $({\cal E}_{1 1}, {\cal E}_{22}, A^+, A^-)$ is a Morita context, so ${\cal E} = {\cal E}_0\oplus {\cal E}_1$ with even part ${\cal E}_0  =  {\cal E}_{11}\oplus{\cal E}_{22}$ and odd part ${\cal E}_1 ={\cal E}_{12}\oplus{\cal E}_{21}$ is a Morita superalgebra according to the definition introduced in \cite[1.4 (III)]{Liesup}).

 The pair $({\cal E},e)$ (or simply, the associative algebra ${\cal E}$ if the idempotent $e$ is understood) is  termed the {\sl standard imbedding} of the associative pair $A$. The {\sl associative  envelope}  of an associative pair $A=(A^+,A^-)$ is the subalgebra ${\cal A}$  of its standard imbedding ${\cal E}$ generated by the odd part of the superalgeba ${\cal E}$.

 The associative envelope ${\cal A}$ of $A$ is  an essential ideal of the standard imbedding ${\cal E}$
 and $A=({\cal A}_{12},{\cal A}_{21})$, where ${\cal A}_{ij}=e_i {\cal A} e_j$, $i,j=1,2$ (Indeed the Peirce projections $\pi_{ij}:{\cal E}\to {\cal E}_{ij}$ can  be restricted to ${\cal A}\to {\cal A}_{ij} = {\cal E}_{ij}\cap {\cal A}$, $i,j=1,2$.)
 \end{apartado}

\begin{remark} Notation and terminology for what we have referred to as the standard imbedding and the associative envelope of associative pairs have been rather interchangeably used in the literature.
A careful review of the different references mentioned in the present paper should allow the reader to tackle this ambiguous usage.  This will be, for instance, the case for the two  references  \cite{fggs-goldie pares} and \cite{ft}, where despite the used notation, the authors deal with the standard imbedding of associative pairs. See  \cite[p.~2998]{fggs-goldie pares}  and \cite[3.3]{ft} for more details. In \cite{gs-LQ AP MP}
standard imbedding and  associative  envelope appear  as introduced in \ref{subsect standard imbedding}. Different notations are used for instance in \cite{CGM-primitive AP} or \cite{pares-PI}. We remark here that since ${\cal A}$ is  an essential ideal of  ${\cal E}$, all results proved along the paper will hold for both   ${\cal E}$ and   ${\cal A}$.\end{remark}

  \begin{apartado}\label{subsect involution}   An {\sl involution} in an associative pair  $A=(A^+,A^-)$ (sometimes named a polarized involution)
  is a pair of  linear mappings $\ast :A^\sigma\to A^\sigma$
such that  $(x^\ast)^\ast=x$ and  $(
x y z)^\ast=  z^\ast y^\ast x^\ast $ for all
$x, z\in A^\sigma$, \  $y\in A^{-\sigma}$, $\sigma =\pm$. Every (polarized) involution of an  associative pair  $A=(A^+,A^-)$  extends uniquely to an involution on its standard imbedding ${\cal E}$   wich coincides with $\ast$ on ${\cal E}_{1 2} = A^+$ and ${\cal E}_{2 1}= A^-$, and satisfies $e_1^\ast=e_2$ \cite[3.2]{ft}.
\end{apartado}

   \begin{apartado}\label{subsect ideal AP} A {\sl left ideal}  of an
  associative pair $A=(A^+,A^-)$ is a
$\Phi$-module  $L$ of $A^\sigma $ such that
$A^\sigma  A^{-\sigma  }L\subseteq L$, $\sigma =\pm$. Right ideals
are defined similarly. A {\sl two-sided ideal} is simultaneously a left and a right ideal.    A pair $I=(I^+,I^-)$ of
two-sided ideals of $A$ is   an {\sl ideal} if
$A^\sigma I^{-\sigma}A^\sigma\subseteq I^\sigma$, $\sigma=\pm$.
  An associative pair  $A=(A^+,A^-)$  is {\sl semiprime} if and only if $I^\sigma A^{-\sigma} I^\sigma=0$, $\sigma=\pm$, implies $I=0$ and {\sl prime} if $I^\sigma A^{-\sigma} J^\sigma=0$, $\sigma=\pm$, implies $I=0$ or $J=0$, for any ideals  $I$, $J$ of $A$. If $A$ is semiprime, then for any ideal $I=(I^+,I^-)$ of $A$ it follows easily that $I^+=0$ if and only if $I^-=0$ (see for instance \cite[p.~2992]{fggs-goldie pares}).  Primeness
implies nondegenerancy ($a^\sigma A^{-\sigma}a^\sigma=0$ implies $a^\sigma=0$, $\sigma=\pm$) and semiprimeness is equivalent to nondegenerancy.
    \end{apartado}

    \begin{apartado}\label{subsect annhilation AP-SI} Let $A$ be an associative pair with standard imbedding ${\cal E}$. Then
    \begin{align*}
    & x_{11}{\cal E}_{12}={\cal E}_{21}x_{11}=0 \quad \Rightarrow  \quad  x_{11}=0,\\
      & x_{22}{\cal E}_{21}={\cal E}_{12}x_{22}=0  \quad\Rightarrow  \quad  x_{22}=0.
    \end{align*}
    If $A$ is semiprime, the above conditions reduce to:
    \begin{align*}
    & x_{11}{\cal E}_{12}=0\quad  \Rightarrow  \quad  x_{11}=0,\\
      & x_{22}{\cal E}_{21}=0\quad  \Rightarrow  \quad  x_{22}=0,
    \end{align*}
    or, equivalently,
     \begin{align*}
    &  {\cal E}_{21}x_{11}=0  \quad\Rightarrow   \quad x_{11}=0,\\
      &  {\cal E}_{12}x_{22}=0 \quad \Rightarrow \quad   x_{22}=0.
    \end{align*}
          \end{apartado}

    \begin{apartado}\label{subsect regularity conditions}    (Semi)primeness of any associative pair is equivalent to that of its standard imbedding \cite[Proposition 4.2]{fggs-goldie pares}. This result stems from the correspondence between ideals of the associative pair  $A=(A^+,A^-)$  and ideals of its standard imbedding ${\cal E}$ (see \cite[Proposition 4.1]{fggs-goldie pares}).
   \end{apartado}

\begin{apartado}\label{subsect annihilator}
 Let  $A=(A^+,A^-)$ be an associative pair. For any subset  $X\subseteq A^\sigma$, $\sigma=\pm$, the the {\sl left} and {\sl right annihilators} of $X$ in $A$ are the sets
 \begin{align*}
 &lann_A(X)=\{b\in A^{-\sigma}\mid  bXA^{-\sigma}=A^\sigma bX=0\},\\
 & rann_A(X)=\{b\in A^{-\sigma}\mid Xb A^\sigma=A^{-\sigma}Xb=0\},
  \end{align*}
   and if $A$ is semiprime, then
  \begin{align*}
 &lann_A(X)=\{b\in A^{-\sigma}\mid   A^\sigma bX=0\},\\
 & rann_A(X)=\{b\in A^{-\sigma}\mid  A^{-\sigma}Xb=0\}.
  \end{align*} which  are left and right ideals of $A$ respectively. 
  The  {\sl annihilator } of $X$ is $ann_A(X)= lann_A(X)\cap rann_A(X)$. If $I=(I^+,I^-)$ is an
   ideal of a semiprime associative pair $A$, then
     the {\sl annihilator}  $ ann_A(I)=(ann_A(I^+), ann_A(I^-))$ of   $I$ is
   $$ ann_A(I^\sigma)=\{x\in A^{-\sigma}\mid xI^\sigma x=0\}. $$  Moreover  $ ann_A(I)$
   is itself an ideal of $A$, and   $I^\sigma \cap ann_A(I^{-\sigma})=0$, $\sigma=\pm$ \cite[Proposition 2.2]{fggs-goldie pares}.
   \end{apartado}

   \begin{apartado}
 An ideal $I$ of an associative pair $A$ is {\sl  essential } if  $I\cap J\neq0$ for any nonzero ideal $J$ of $A$. It follows from \cite[Proposition 2.2]{fggs-goldie pares} that essential ideals of semiprime associative pairs are exactly those ideals of $A$ whose annihilator vanishes.     \end{apartado}

   \begin{apartado}\label{subsect socle AP}    The {\sl socle} of a semiprime associative pair $A $  consists of the pair of subsets $Soc(A)=(Soc(A^+), Soc(A^-))$, where $Soc(A^\sigma)$, $\sigma=\pm$,  is the sum of all minimal right ideals of $A$.   $Soc(A)$ is an ideal of $A$, and if $Soc(A)\neq0$, then it is a direct sum of simple ideals.   If $A$ is prime, $Soc(A)$ is a simple ideal contained in every nonzero ideal of $A$ \cite[Theorem 1]{CGCM-sp AP miin-i}. Elements of the socle of semiprime associative algebras and pairs are  von Neumann regular \cite[Theorem 1]{loos-jp}.
We refer to  \cite{CGM-primitive AP}, \cite{fg-socle ATS}  and \cite{ft} for descriptions of prime associative pairs with nonzero socle.   Prime associative pairs with involution having nonzero socle together with their involutions  were described in \cite[Theorem 3.14]{ft}. \end{apartado}

      \begin{apartado}\label{subsect socle SI}   There is a good relation between the socle of a semiprime associative pair $A$ and that of its standard imbedding ${\cal E}$:
    $Soc(A^+)=Soc({\cal E})\cap A^+=e_1Soc({\cal E})e_2$,
   and similarly $Soc(A^-)=Soc({\cal E})\cap A^-=e_2Soc({\cal E})e_1$ \cite[Proposition 3.4(4)]{ft}. The standard imbedding of $Soc(A)$ can be identified with the ideal $Soc({\cal E})$ of the standard imbedding ${\cal E}$ of $A$.
   \end{apartado}

 \begin{apartado}\label{subsect primitivity} A  pair of $\Phi$-modules $M=(M^+,M^-)$ is a {\sl right  $A$-module},  for an associative pair $A$, if $M$ is endowed with a pair of  $\Phi$-bilinear maps
\[\begin{array}{ccc}   M^\sigma \times A^{-\sigma} & \to & M^{-\sigma}\\
\quad (m,x) & \mapsto & mx\end{array}\]
satisfying $((mx)y)z=m(xyz)$ for all $m\in M^\sigma$, $x,z\in A^{-\sigma}$ $y\in A^\sigma$, $\sigma=\pm$. Left $A$-modules are defined similarly. A right $A$-module $M=(M^+,M^-)$ is {\sl irreducible} if $M^{-\sigma}A^\sigma\neq0$, $\sigma=\pm$, and it   contains no proper submodules (different from 0 and $M$ itself) and {\sl  faithful} if $M^{-\sigma}x=0$ implies $x=0$ for any $x\in A^\sigma$, $\sigma=\pm$.
 An associative pair  $A=(A^+,A^-)$  is {\sl right primitive} if it has a  faithful irreducible right $A$-module. In  \cite[Theorem 1]{CGM-primitive AP}, it is proved that a Density Theorem holds for primitive associative pairs, and that an associative pair is primitive if and only if so is its standard imbedding.
 \end{apartado}

   \begin{apartado}\label{subsect AP regularity}   Primitive associative pairs are prime, and associative pairs with nonzero socle are primitive if and only if they are prime  \cite[2.8]{CGM-primitive AP}.
A Structure Theorem for primitive associative pairs with nonzero socle was given in \cite[Theorem 2]{CGM-primitive AP}.  See also \cite[Theorem 3.9]{ft}.
\end{apartado}

 \begin{apartado}\label{subsect def local algebra}
The {\sl local algebra} of an associative pair $A$ at an element $a\in A^{-\sigma}$ is the quotient algebra $A^\sigma_a=(A^\sigma)^{(a)}/Ker \, a $ of the {\sl $a$-homotope algebra} $(A^\sigma)^{(a)}$ of $A$ (the associative algebra over the $\Phi$-module $A^\sigma$ with product $x\cdot y=xay$ for all $x,y \in A^\sigma$) over the ideal $Ker\, a=\{x\in A^\sigma\mid axa=0 \}$ of $(A^\sigma)^{(a)}$ \cite{meyberg,pi-i}. Local algebras of associative pairs interact well with standard imbeddings: $A_a\cong {\cal E}_a$ for all $a\in A^{-\sigma}$ \cite[Proposition 3.4(3)]{ft}.

As for regularity conditions and their interaction to local algebras, we recall from \cite[Proposition 5.2]{fggs-goldie pares} the following facts: Local algebras of semiprime associative pairs are semiprime associative algebras, an  associative pair $A$ is prime   if and only if all its local algebras at nonzero elements are prime, and if $A$ is simple, then so are all its local algebras at nonzero elements \cite[Proposition 5.2]{fggs-goldie pares}.
 \end{apartado}

\begin{apartado}
Associative pairs satisfying polynomial identities were studied in \cite{pares-PI}.
We denote by $FAP(X)$ the {\sl free associative pair} over $\Phi$ on  indeterminates $X=(X^+,X^-)$, which is the subpair of the pair $(FA(X^+\cup
X^-),FA(X^+\cup X^-))$ obtained by doubling the free  associative
algebra $FA(X^+\cup X^-)$, generated by $(X^+,X^-)$.
The universal property of $FAP(X)$ makes possible to evaluate any pair polynomial $f_\sigma(x_1^+,\ldots x_n^+,x_1^-,\ldots,x_n^-)$ on an associative pair $A$, by assigning fixed values $x_i^\sigma=a_i^\sigma\in A^\sigma$.
  An associative polynomial $f_\sigma\in FAP(X)^\sigma$ is a {\sl polynomial identity} of an associative pair $A$, if $f_\sigma$ is monic (i.e.  some of its leading monomials has coefficient 1) and all evaluations of $f_\sigma $ on $A$ vanish.
   Similarly we can consider  {\sl $\ast$-polynomials}   $ p_\sigma(x_1^+,\ldots, x_n^+, (x^+_1)^\ast,\ldots, (x^+_n)^\ast,
   x_1^-,\ldots, x_n^-, (x^-_1)^\ast,\ldots, (x^-_n)^\ast)$  and {\sl $\ast$-polynomial identities}.  We will  say that an associative pair is $PI$ (or that it is an {\sl associative PI-pair}) if it satisfies a polynomial identity, and similarly one defines  {\sl $\ast$-PI associative pairs}.
\end{apartado}

For a primitive pair, satisfying a a ($\ast$-)polynomial identity ensures the existence of nonzero socle:

 \begin{proposicion} {\rm \cite[Proposition 3.4, Theorem 3.6]{pares-PI}} Let $A$ be a primitive associative pair.
   \begin{enumerate}\item[(i)] If $A$ is PI, then $A$ has nonzero socle.
      \item[(ii)] If $A$ has an involution $\ast$, and is $\ast$-PI, then $A$ has nonzero socle.
   \end{enumerate} Moreover in either case, $A$ is simple and has finite capacity.
\end{proposicion}

\begin{remark} The capacity of PI (or $\ast$-PI) primitive pairs is bounded by a constant depending only on the degree of the polynomial identity (see \cite[Theorem 3.6]{pares-PI}).
\end{remark}

The following analogue of Amitsur's theorem for associative pairs  with involution was also proved in \cite{pares-PI}.

 \begin{teorema} {\rm \cite[Theorem 3.9]{pares-PI}} Let $A$ be an associative pair with involution $\ast$. If $A$ has a $\ast$-polynomial identity of degree $m$,   there exists a positive integer $k$ such that every local algebra of $A$ satisfies the polynomial identity $S^k_{2m}$. Moreover, if $A$ is semiprime,  every local algebra satisfies the standard identity $S_{2m}$.
\end{teorema}

 \begin{apartado} The notion of   PI-element for associative pairs  was introduced in  \cite{pi-i}: An element $a \in A^{-\sigma}$ of an associative pair $A$ is a {\sl PI-element} if the local algebra $A^\sigma_a$ satisfies a polynomial identity. Then, the pair   $PI(A)=(PI(A^+),PI(A^-))$,  where
 $PI(A^\sigma)$ denotes the set of all PI-elements of $A^\sigma$, is an ideal of the associative pair $A$ \cite[Proposition 1.6]{pi-i}.
 \end{apartado}

\begin{proposicion}\label{prop PI-ideal}  Let $A$ be a semiprime associative pair. Then  $PI(A)=PI({\cal E})\cap A$.  \end{proposicion}
\begin{proof}
This follows from the relation equality  \ref{subsect def local algebra}  between the local algebras of the standard imbedding ${\cal E}$ of $A$  at elements of the pair $A$, and the local algebras of the associative pair $A$.
\end{proof}

 \begin{apartado} \label{subsect centroid AP}   The
centroid $\Gamma (A)$ of an associative pair $A=(A^+,A^-)$ is the set of all pairs
$T=(T^+,T^-)\in End_\Phi(A^+)\times End_\Phi(A^-)$  satisfying:
$$T^\sigma(x^\sigma y^{-\sigma }z^\sigma )=
T^\sigma (x^\sigma) y^{-\sigma }z^\sigma =x^\sigma T^{-\sigma
}(y^{-\sigma })z^\sigma=x^\sigma y^{-\sigma }T^\sigma (z^\sigma)$$
for all $x^\sigma, z^\sigma\in A^\sigma $, $y^{-\sigma }\in
A^{-\sigma }$, $\sigma=\pm$.
The centroid of a  semiprime associative pair  is a
commutative ring. If $A$ is prime, then $\Gamma (A)$ is  a domain acting faithfully on $A$,   and it is a field if $A$ is simple.
 \end{apartado}

\begin{apartado} The extended centroid and the central closure of associative rings were introduced by Martindale for prime rings \cite{martindale prime rings} and generalized to semiprime rings by  Amitsur \cite{amitsur sprime rings}.   The nonassociative case was considered in
  \cite{EMO-centroid-prime} and \cite{BM-centroid-sprime}, and
   associative rings with involution  were dealt with in \cite{BM-centroid-prime involution} and \cite{BM-centroid-sprime involution}.
  In \cite{pi-ii} the notions of extended centroid and the analogue to
  the central closure, called there extended central closure (since in that context there was already a notion of central closure) were introduced for
   Jordan systems (algebras, pairs and triple systems) on arbitrary rings of scalars.
\end{apartado}

 \begin{apartado}\label{subsect extended centroid rings}  We briefly recall now the definition of the extended centroid of an associative algebra.

 Let $I$ be an ideal of an associative algebra $R$.  If $f:I\to R$ is a homomorphism of $R$-bimodules, and $I$ is an essential ideal of $R$, then $f$ will be called a {\sl permissible map}. We will write it as the pair $(f,I)$, since we will make use of the restrictions of $f$ to smaller essential ideals, so it is convenient to have the domain explicitly displayed.

 The extended centroid of a semiprime ring $R$ is the direct limit $\displaystyle {\cal C}(R) = \lim_{\rightarrow} Hom_R(I,R)$ over the filter of essential ideals of $R$ with the operations naturally inherited from $R$. Explicitly, as a set, ${\cal C}(R)$ consists of the equivalence clases of permissible maps $\overline{(f,I)}$ under the equivalence relation $(f,I)\thicksim (g,J)\quad\hbox{if $f_{\mid L}=g_{\mid L}$ for some essential ideal $L\subseteq I\cap J$} $.

 The operations in ${\cal C}(R)$ are defined by
 $ \displaystyle \overline{(f,I)}+ \overline{(g,J)}=  \overline{(f+g,I\cap J)},$
 and
 $\displaystyle \overline{(f,I)} \, \cdot\, \overline{(g,J)}=  \overline{(fg ,g^{-1}(I))}$. As a result, ${\cal C}(R)$ becomes a commutative von Neumann regular unital ring \cite[Theorem 2.5]{BM-centroid-sprime} which is called the {\sl extended centroid} of the associative algebra $R$. The corresponding scalar extension ${\cal C}(R)R$, called the {\sl central closure} of $R$, remains semiprime, is generated as a  ${\cal C}(R)$-module by $R$, and  is centrally closed \cite[Theorem 2.15]{BM-centroid-sprime}. We also recall here that the extended centroid of a semiprime ring $R$ is the center of its maximal  right (and   left) ring of quotients (also that of its symmetric ring of quotients) \cite[Remark 2.3.1]{BMM-rings GPI}. We also refer the reader to \cite[Section~3]{BM- jordan homomorph} for   further information on the construction of the central closure of a semiprime ring.
  \end{apartado}

\section{The extended centroid of semiprime associative pairs}

In this section we first extend the construction of the extended centroid for semiprime  algebras \cite{amitsur sprime rings,BM-centroid-sprime,martindale prime rings}  to semiprime associative pairs, and then relate the extended centroid of a semiprime associative pair to  that of its standard imbedding.

 \begin{apartado} Let $A=(A^+,A^-)$  be an associative pair and let $I=(I^+,I^-)$ be an ideal of $A$. Then $f:I\to A$ is an {\sl $A$-homomorphism} if $f=(f^+,f^-)$ consists of a pair of $\Phi$- linear maps $f^\sigma: I^\sigma\to A^\sigma$, $\sigma=\pm$, satisfying:
\begin{align*}
&f^\sigma(x^\sigma y^{-\sigma} z^\sigma)=  x^\sigma f^{-\sigma}(y^{-\sigma}) z^\sigma,\\
&f^\sigma(y^\sigma x^{-\sigma} z^\sigma)=   f^{ \sigma}(y^{ \sigma})x^{-\sigma}  z^\sigma,\\
&f^\sigma(x^\sigma z^{-\sigma} y^\sigma)=  x^\sigma  z^{-\sigma}   f^\sigma(y^\sigma),
\end{align*}
 for all $x^\sigma, z^\sigma\in A^\sigma$, $y^\sigma\in I^\sigma$, $\sigma=\pm$. \end{apartado}

  \begin{apartado} Note that a pair of $\Phi$- linear maps
 $f=(f^+,f^-)$ is an $A$-homomorphism if and only if it commutes  with all left, right and middle multiplication operators defined by elements of $A$. We denote by $Hom_A(I,A)$,  where $ Hom_A(I,A) =(Hom_A(I^+,A^+),Hom_A(I^-,A^-) )$, the set of all $A$-homo\-morphisms from $I$ to $A$. A pair $(f,I)$ where $f\in Hom_A(I,A)$ will be called a {\sl permissible map} if $I$ is an essential ideal of the associative pair $A$.
 \end{apartado}

\begin{teorema}\label{th extended centroid AP} Let $(f,I)$ and $(g,J)$ be permissible maps of a semiprime associative pair $A$. Then:
 $$(f,I) \thicksim (g,J) \ \hbox{if}\ f_{\mid K}=g_{\mid K} \ \hbox{for some essential ideal} \ K\subseteq I\cap J,$$
defines an equivalence relation in the set of all $A$-permissible maps of $A$.
Then the quotient set ${\cal C}(A)$,   with the operations:
$$ \overline{(f,I)}+ \overline{(g,J)}=  \overline{(f+g,I\cap J)},$$
 $$  \overline{(f,I)} \cdot \overline{(g,J)}=  \overline{(fg ,g^{-1}(I))},$$
is a commutative, von Neumann regular unital ring.
\end{teorema}

\begin{proof}
This is straightforward, arguing as for the corresponding results on algebras  \cite{BM-centroid-sprime}, mentioned before.
\end{proof}

 \begin{apartado}  We will refer to ${\cal C}(A)$  as the {\sl  extended centroid} of the semiprime associative pair $A$. Clearly ${\cal C}(A)$  contains a copy of the centroid $\Gamma(A)$ of $A$.
 \end{apartado}

 \begin{apartado} Our aim next is  to relate  the extended centroid  of  a  semiprime associative pair   to that of its standard imbedding. Since we will be simultaneously dealing with ideals of associative pairs and of their standard imbeddings, we will denote by $I, J, K,\ldots $ the associative pair ideals and by  $ {\cal I} , {\cal J}, {\cal K}, \ldots$ the algebra ideals.
 \end{apartado}

\begin{remark} It is obvious that if $I$ is an essential ideal of an associative algebra $R$,   the extended centroid ${\cal C}(I)$ can be identified with the extended centroid ${\cal C}(R)$ by considering the homomorphism induced by the restriction of permissible maps $(f, L)\mapsto (f_{\vert I\cap L}, L\cap I, )$. Therefore it does not matter  whether we work with the standard imbedding or the associative envelope of associative pairs.
\end{remark}

\begin{lema}\label{lemma ideals 1} Let $ {\cal E}$ be  the standard imbedding of  an associative pair $A=(A^+,A^-)$.
\begin{enumerate}
\item[(i)] If ${\cal I}$  is a nonzero ideal of  $ {\cal E}$, then $I=({\cal I}\cap A^+ ,{\cal I}\cap A^- ) =({\cal I}_{12},{\cal I}_{21})$  is a nonzero ideal of $A$.
\item[(ii)]   If $I=(I^+,I^-)$ is a nonzero ideal of $A$, then, the ideal of ${\cal E}$ generated by $I$ is
$${\cal I}=\big(I^+A^-+A^+I^-\big)\oplus I^+ \oplus I^- \oplus\big(I^-A^++A^-I^+\big)$$  which
is nonzero. Moreover if $A$ is semiprime and $I$ is essential in $A$, then  so is $ {\cal I}$ in $ {\cal E}$.
\end{enumerate}
\end{lema}

  \begin{proof} For (i)  see
  \cite[Proposition 4.1(i)]{fggs-goldie pares}.

  (ii)    This is straightforward from (i) and \ref{subsect standard imbedding}.
    \end{proof}

  \begin{lema}\label{lemma ideals 2} Let $ {\cal E}$ be the  standard imbedding of  a semiprime associative pair $A=(A^+,A^-)$.
  If ${\cal I} $  is an essential ideal of  $ {\cal E}$, then $I=({\cal I}\cap A^+ ,{\cal I}\cap A^- ) =({\cal I}_{12},{\cal I}_{21})$ is essential in $A$.
    \end{lema}

   \begin{proof} This easily follows from the previous Lemma. \end{proof}

 \begin{lema}\label{lema Pierce} Let
 $A=(A^+,A^-)$
be  a semiprime associative pair. Then for any
 $ {\cal E}$-homomorphism $(f,{\cal I}  )$    of its standard imbedding  $ {\cal E}$, we have
   $f  \pi_{ij}=\pi_{ij}  f$, for all $i,j=1,2$ (i.e.    $ {\cal E}$-homomorphisms are compatible with the Peirce decomposition of  $ {\cal E}$).
\end{lema}

\begin{proof} Write  ${\cal I}={\cal I}_{11}\oplus {\cal
I}_{12}\oplus{\cal I}_{21}\oplus {\cal I}_{22}$, where  ${\cal
I}_{ij}=e_i{\cal I}e_j={\cal I}\cap {\cal E}_{ij}$, and take
$x=x_{11}+x_{12}+x_{21}+x_{22}\in {\cal I}$.  Then, since $x_{ij}\in {\cal
I}_{ij}\subseteq {\cal
I}$,  we have
$$f\pi_{ij}(x)=f(x_{ij})=f(e_ix_{ij}e_j)=f(e_ixe_j)=e_if(x)e_j=\pi_{ij}f(x). $$
Hence
 $f\pi_{ij}(x)=\pi_{ij}f(x)$   for all $i,\ j\in \{1,2\}$, $x\in {\cal
I}$.
\end{proof}

\begin{lema}\label{lema permissibles restric}  Let  $A$ be a semiprime associative pair, and let
  $\lambda_1 =\overline{(f,{\cal I})}$ and $
\lambda_2=\overline{(g,{\cal J})}$ be  elements of the extended centroid ${\cal C(E)}$ of the standard imbedding  ${\cal E}$ of $A$.
If $f_{\mid {\cal I}\cap{\cal J}\cap
A}=g_{\mid{\cal I}\cap{\cal J}\cap A}$, then $\lambda_1
=\lambda_2 $ in  ${\cal C(E)}$.
\end{lema}

\begin{proof} We note first that replacing ${\cal I}$  and ${\cal J}$ by ${\cal I}\cap  {\cal J}$, we can assume ${\cal I}={\cal J}$ \cite[Corollary 2.3]{BM-centroid-sprime}. Write then  $\lambda_1 =\overline{(f,{\cal I})}$ and $
\lambda_2=\overline{(g,{\cal I})}$ where ${\cal I}={\cal I}_{11} \oplus {\cal I}_{12}\oplus {\cal I}_{21}\oplus {\cal I}_{22}$ is an essential ideal of ${\cal E}$ and assume $f(y_{12})=g(y_{12}) $ and $f(y_{21})=g(y_{21}) $  for all $y_{12}\in {\cal I}_{12}$ and all $y_{21}\in {\cal I}_{21}$.

Take now $y_{11}\in {\cal I}_{11}$  and $a_{12}\in {\cal E}_{12}=A^+$. Then since $y_{11}a_{12}\in  {\cal I}_{11}{\cal E}_{12}\subseteq {\cal I}_{12}={\cal I}\cap A$, we have
$$\big(f(y_{11})-g(y_{11})\big)a_{12}=f(y_{11})a_{12}-g(y_{11})a_{12}=f(y_{11}a_{12})-g(y_{11}a_{12})=0.$$
Hence
 $\big(f(y_{11})-g(y_{11})\big)A^+=0$, where by   Lemma \ref{lema Pierce},   $f(y_{11}), g(y_{11})\in  {\cal E}_{11}$. Therefore  by the semiprimeness of $A$ (see \ref{subsect annhilation AP-SI}),  this implies $f(y_{11})=g(y_{11})$ for all $y_{11}\in {\cal I}_{11}$. Thus we have
 $f_{\mid   {\cal I}_{11}}=g_{\mid{\cal
I}_{11} }$ and similarly one proves that $f_{\mid   {\cal I}_{22}}=g_{\mid{\cal
I}_{22} }$. Hence
   $\lambda_1=\lambda_2  $.
\end{proof}

\begin{corolario}\label{corollary permissibles restric}
 Let  $A$ be a semiprime associative pair and let
  $\lambda =\overline{(f,{\cal I})}  $ be  an element of the   extended centroid  ${\cal C(E)}$ of the standard imbedding  ${\cal E}$ of $A$.
If   $f_{\mid {\cal I}\cap A}=0$ then $\lambda=0$.
\end{corolario}

\begin{proof} This result  is a particular case of Lemma \ref{lema permissibles restric}.\end{proof}

\begin{teorema}\label{th iso centroids} Let $A$ be a semiprime associative pair with standard imbedding ${\cal E}$. Then
 ${\cal C}(A)\cong {\cal C}({\cal E})$.
\end{teorema}

\begin{proof}
Take first an element $\lambda =\overline{(f,{\cal I})}\in {\cal C}({\cal E})$, where ${\cal I }$ is an essential ideal of ${\cal E}$ and $f:{\cal I }\to {\cal E}$
is a permissible map. Write $I=(I^+,I^-)=({\cal I}_{12},{\cal I}_{21})={\cal I}\cap A$, which is an essential ideal of the associative pair  by Lemma \ref{lemma ideals 2}, and consider $g=(g^+,g^-)$ where $g^+=f_{|{\cal I}_{12}}$ and $g^-=f_{|{\cal I}_{21}}$.
By Lemma \ref{lema Pierce},  $g=(g^+,g^-)$ consists of  pair of linear maps with $g^\sigma: I^\sigma \to A^\sigma$, $\sigma=\pm$.

We note that from the previous lemmas it follows that the mapping
 $\Psi: {\cal C}({\cal E})\to {\cal C}(A)$ given by
 $\Psi(\lambda )=\overline{(g,I)}$ is a well-defined ring homomorphism.   Besides, by Corollary \ref{corollary permissibles restric}, $\Psi$ is easily seen to be injective.

 Conversely, take now $\mu =\overline{(g,I)}\in {\cal C}(A)$. By Lemma \ref{lemma ideals 1}:
$\displaystyle {\cal I}=\big(I^+A^-+A^+I^-\big)\oplus I^+\oplus I^-\oplus
\big(I^-A^++A^-I^+\big)$ is an essential ideal of ${\cal E}$ which satisfies  ${\cal I}\cap A=I$.  We now define a linear map $f: {\cal I}\to {\cal E}$ as follows:
  \begin{enumerate}
    \item[(i)] $\sum_{i=1}^n g^+(y_i^+)a_i^-+\sum_{j=1}^m b_j^+g^-(x_j^-)$ for all  $\sum_{i=1}^n y_i^+ a_i^-+\sum_{j=1}^m b_j^+ x_j^- \in{\cal I}_{11}= I^+A^-+A^+I^-$, $y_i^\sigma, x_j^\sigma\in I^\sigma$, $a_i^\sigma, b_j^\sigma\in A^\sigma$,  $\sigma=\pm$,
    \item[(ii)] $g^\sigma(y^\sigma)$ for all $y^\sigma\in I^\sigma$, $\sigma=\pm$,
    \item[(iii)] $\sum_{k=1}^p g^-(z_k^-)c_k^++\sum_{l=1}^q  d_l^-g^+(t_l^+)$ for all  $\sum_{k=1}^p  z_k^- c_k^++\sum_{l=1}^q d_l^- t_l^+  \in{\cal I}_{22}= I^-A^++A^-I^+$, $z_k^\sigma, t_l^\sigma\in I^\sigma$, $c_k^\sigma, d_l^\sigma\in A^\sigma$,  $\sigma=\pm$,
  \end{enumerate}
  and extended it to ${\cal I}$ by linearity.
  Again, the mapping
   $\Phi:\ {\cal C}(A)\to {\cal C}({\cal E})$  given by
    $\Phi(\mu )=\overline{(f,{\cal I})}$ is easily seen to be a well-defined injective ring homomorphism.

      Finally, the Theorem follows by noticing that $\Psi$ and $\Phi$ are mutually inverses.
\end{proof}

Even though the direct proof (in a pair environment) of the following result is is also possible, we obtain it as consequence of Theorem \ref{th iso centroids}.

\begin{corolario} The extended centroid ${\cal C}(A) $ of a prime associative pair $A$ is a field.\end{corolario}

\begin{proof} It follows from \cite[Theorem 2.1]{EMO-centroid-prime} since, as noted in
\ref{subsect regularity conditions}, an associative pair $A$ is prime if and only if its standard imbedding is a prime associative algebra.
\end{proof}

\section{Closure of semiprime associative pairs}

The extended centroid of a  semiprime associative pair gives rise to a scalar extension that will be called the central closure. This section sketches the construction of the central closure for semiprime associative pairs. For more explicit details the reader is referred to \cite{BM-centroid-sprime,BM- jordan homomorph}. As expected, central closures and standard imbeddings will be commuting constructions for semiprime associative pairs.

\begin{apartado}\label{subsect central closure AP} We define the {\sl central
 closure} ${\cal C}(A)A$ of a semiprime associative pair $A$ to be the quotient pair of the free scalar extension
$$ {\cal C}(A)\otimes_\Phi A=({\cal C}(A)\otimes_\Phi A^+,{\cal C}(A) \otimes_\Phi A^-)$$ by the pair ideal $R=(R^+,R^-)$, being $R^\sigma$ the linear span of all elements of the form $ \mu   \otimes y^\sigma -1\otimes g^\sigma (y^\sigma)$, where $ \mu\in {\cal C}(A)$, with
 $\mu =\overline{(g,I)}\in  {\cal C}(A)$, $g=(g^+,g^-)$ and $y^\sigma\in I^\sigma$, $\sigma=\pm$. Then ${\cal C}(A)A=(  {\cal C}(A)A^+,{\cal C}(A)A^- )$ and elements of ${\cal C}(A)A^\sigma$ will be written as $a^\sigma =\sum_{i=1}^n \lambda_ix_i^\sigma$, with $\lambda_i=\overline{(g_i, I_i)}$, being $g_i=(g^+_i,g^-_i)$ a pair of  $A$-homomorphisms,  $I_i=(I^+_i,I^-_i)$ an  essential ideal of $A$ and    $x_i^\sigma\in A^\sigma$,  $i=1,\ldots,n$, $\sigma=\pm$.
\end{apartado}

\begin{teorema}\label{th central closure AP} The central   closure  ${\cal C}(A)A$ of a semiprime associative pair $A$ is a tight scalar extension of $A$, and therefore it is  a semiprime associative pair. Moreover, if $A$ is prime, so is ${\cal C}(A)A$.
\end{teorema}

\begin{proof} This follows as the corresponding algebra result of  \cite{BM-centroid-sprime}, with the obvious changes for associative pairs.
\end{proof}

\begin{apartado} If $A$ is an associative pair with standard imbedding ${\cal E}$, then the pair of orthogonal idempotents $e_1$ and $e_2$, with  $e_1+e_2=1$, and such that $A\cong ({\cal E}_{12},{\cal E}_{21})$,  also induces a Peirce decomposition on the central
 closure ${\cal C(E)E}$ of the associative algebra ${\cal E}$. We next prove that the  associative pair given by the off-diagonal Peirce components of ${\cal C(E)E}$ corresponds to the central closure  ${\cal C}(A)A$ of the associative pair $A$.\end{apartado}

\begin{teorema}\label{th iso  central closures} Let $A$ be a semiprime associative pair with standard imbedding ${\cal E}$. Then the standard imbedding of the central closure  ${\cal C}(A)A$ of $A$ is isomorphic to the central   closure   ${\cal C(E)E}$ of the standard imbedding ${\cal E}$ of $A$.
\end{teorema}

\begin{proof} It suffices to consider the induced Peirce decomposition of  ${\cal C(E)E}$  given by the pair of orthogonal idempotents $e_1$ and $e_2$ of $ {\cal E}$, taking  into account \ref{subsect standard imbedding} and Theorem \ref{th iso centroids}.
\end{proof}

\begin{corolario}  The central closure   ${\cal C}(A)A$ of  a semiprime associative pair $A$ is closed over  ${\cal C}(A)$.  \end{corolario}
\begin{proof}
It follows from Theorem \ref{th iso central closures}, Theorem \ref{th iso centroids} and \cite[Theorem~2.15(c)]{BM-centroid-sprime}.
\end{proof}

\begin{remark}\label{remark Martindale AQ}
We note here that ${\cal C}(A)A= (e_1{\cal C(E)E} e_2, e_2 {\cal C(E)E} e_1 )$.  This representation of the central closure ${\cal C}(A)A$ of  a semiprime associative pair $A$ as the off-diagonal Peirce spaces of the central closure of its standard imbedding  ${\cal E}$ agrees with the
  definition given in \cite[Definition 2.11]{gs-LQ AP MP} of the maximal left quotient pair  of   associative pairs without total zero divisors (and, in particular, for semiprime associative pairs). The same approach can be indeed applied to  other constructions of pairs of quotients of a  semiprime associative pair $A$, as for instance, the maximal pair of symmetric quotients $Q_\sigma(A)$   or the Martindale symmetric ring of quotients $Q_s(A)$ of $A$
\cite{pares-quot}.
\end{remark}

\section{Associative pairs with involution}

In the present section we review the versions of the  results on extended centroids and central  closures for  semiprime  pairs  with a (polarized) involution $\ast$. Recall (see \ref{subsect involution}) that such an involution extends uniquely to an involution also denoted by $\ast$ on the standard imbedding ${\cal E}$ of $A$, with $e_1^\ast=e_2$.

\begin{apartado}\label{subsect extended centroid rings inv}  Involutions of a semiprime ring $R$ extend easily to its extended centroid. Indeed, given $\lambda=\overline{(f,I)}\in {\cal C}(R)$, where $f$ is an $R$-homomorphism and $I$ an essential
$\ast$-ideal of $R$, it suffices to define
  $\lambda^\ast=\overline{(f^\ast,I  )}$ where $f^\ast(y)=(f(y^\ast))^\ast$ for all $y\in I$  \cite[p.~1125]{BM-centroid-sprime}.  We recall here that for any semiprime ring with involution, the filter of essential ideals is equivalent to the filter of essential $\ast$-ideals. Then for a semiprime ring $R$ the {\sl $\ast$-extended centroid }
   ${\cal C}_\ast(R)$ of $R$, defined as the set of
  all symmetric elements of ${\cal C}(R)$, is a unital commutative ring.  Similarly, involutions of $R$ also extend  to the central closure ${\cal C}(R)R$ of $R$  \cite{BM-centroid-prime involution, BM-centroid-sprime involution}.
  \end{apartado}

 \begin{apartado} Although  the subring of fixed elements of the extended centroid of a semiprime ring with involution was already considered in \cite{BM-centroid-sprime}, the
   notion of  {\sl $\ast$-extended centroid} was introduced in \cite{BM-centroid-prime involution} for $\ast$-prime rings, and extended to semiprime rings with involution in  \cite[p. 952]{BM-centroid-sprime involution}, based on the set of equivalence classes of  $\ast$-permissible maps defined on essential $\ast$-ideals (these are permissible maps $f:I\rightarrow R$ defined on an essential $\ast$-ideal: $I^\ast =I$, which commute with the involution: $f(y^\ast) = f(y)^\ast$).
     \end{apartado}

\begin{apartado}\label{subsect central  closure rings inv}  The {\sl $\ast$-extended centroid }
   ${\cal C}_\ast(R)$ of a semiprime ring with involution gives rise to a   scalar extension ${\cal C}_\ast(R)R$, called the {\sl $\ast$-central closure}
  of $R$.  Again, ${\cal C}_\ast(R)R$ is endowed  with an involution $\ast$ defined by  $(\sum_{i=1}^n \lambda_i r_i)^\ast= \sum_{i=1}^n \lambda_i r_i^\ast$, for all   $ \lambda_i\in {\cal C}_\ast(R)$   and $r_i\in R$, $i=1,\ldots,n$.
  If $R$ is $\ast$-prime,   ${\cal C}_\ast(R)$ is a field and  ${\cal C}_\ast(R)R$ is a $\ast$-prime algebra over  ${\cal C}_\ast(R)$   generated by $R$ over ${\cal C}_\ast(R)$. Moreover   ${\cal C}_\ast(  {\cal C}_\ast(R)R)= {\cal C}_\ast(  R)$, i.e. ${\cal C}_\ast(R)R$ is $\ast$-closed \cite[Theorem 4]{BM-centroid-prime involution}. (This can also be obtained through the   symmetric ring of quotients $Q_s(R)$ of $R$ \cite[2.3]{BMM-rings GPI}.)   \end{apartado}

\begin{proposicion}\label{proposition * extended centroid} Let $A$ be a semiprime associative pair with involution $\ast$ and let
$(f,I)$ be  a permissible map of $A$.  Then
 $(f^\ast,I^\ast )$ given by $(f^\sigma)^\ast(y^\sigma)= (f^\sigma((y^\sigma)^\ast ) )^\ast$, for all $y^\sigma\in I^\sigma$,  $\sigma=\pm$,  is permissible and this defines an involution on the extended centroid ${\cal C}(A)$ of $A$.  \end{proposicion}

\begin{proof}  Take  two  permissible  maps $(f,I)$ and $(g,J)$ of $A$, defined on  essential   ideals $I$ and $J$ of $A$ and assume that
 $(f,I)\thicksim(g,J)$. Then  it is straightforward that
    $(f^\ast,I^\ast )\thicksim(g^\ast,J^\ast )$. Moreover $(f^\ast)^\ast=f$. Hence,    ${\cal C}(A)$  being a commutative ring, this defines an involution on ${\cal C}(A)$.
\end{proof}

\begin{teorema}\label{th ast-extended centroid AP} Let $A=(A^+,A^-)$ be a    semiprime associative pair  with an involution $\ast$. The set
$ {\cal C}_\ast(A)$   of all symmetric elements of the extended centroid  $ {\cal C} (A)$ of $A$ with respect to   the involution defined in Proposition \ref{proposition * extended centroid}  forms a commutative unital ring. Moreover $ {\cal C}_\ast(A)$ is a field if $A$ is $\ast$-prime.
\end{teorema}

\begin{proof}   $ {\cal C}_\ast(A)$ is a commutative unital ring as a result of Theorem \ref{th extended centroid AP} and Proposition \ref{proposition * extended centroid}.  Note also that, since for any  $\lambda=\overline{(f, I)}\in {\cal C}_\ast(A)$, both $Ker f=(Ker f^+, Ker f^-)$ and $Im f=(Im f^+, Im f^-)$ are $\ast$-ideals of $A$, if $A$ is $\ast$-prime and $\lambda \neq 0$, then $Ker f$ vanishes. Then $ \mu=\overline{(g, f(I))}$ given by $g^\sigma(f^\sigma (y^\sigma))=y^\sigma$, $\sigma =\pm$,  is an inverse of $\lambda=\overline{(f, I)}$ in ${\cal C}_\ast(A)$ (see \cite[p.~860]{BM-centroid-prime involution}).
\end{proof}

\begin{apartado}    We will refer to   $ {\cal C}_\ast(A)$ as the {\sl $\ast$-extended centroid} of the semiprime associative pair  with involution $A$.
\end{apartado}

\begin{remark}  Given an essential $\ast$-ideal $I=(I^+,I^-)$ of a semiprime associative pair $A$ with involution $\ast$, as for algebras, we say that a pair of $A$-homomorphisms $f=(f^+,f^-)\in
Hom_A(I,A)$ is {\sl $\ast$-permissible} if  $(f^\sigma)^\ast(y^\sigma)=f^\sigma(y^\sigma ) $, for all $y^\sigma\in I^\sigma$,  $\sigma=\pm$. The $\ast$-extended centroid $ {\cal C}_\ast(A)$  of $A$ can be  therefore characterized as the set of all $\ast$-permissible maps in the  extended centroid  $ {\cal C} (A)$  of $A$.
\end{remark}

\begin{teorema}\label{th iso ast-centroid}
  Let $A$ be a semiprime  associative pair   with an involution $\ast$ and standard imbedding ${\cal E}$.  Then   ${\cal C}_\ast(A) \cong {\cal C_\ast(E)}$.
\end{teorema}

\begin{proof}
 We first note that a direct check yields that the maps $\Psi: {\cal C}({\cal E})\to {\cal C}(A)$
   and  $\Phi: {\cal C}({\cal E})\to {\cal C}(A)$  given in Theorem \ref{th iso centroids} are ring  $\ast$-homomorphisms, hence their restrictions define
reciprocal isomorphisms between the $\ast$-extended centroids
    ${\cal C}_\ast(A)$ and    $ {\cal C}_\ast({\cal E})$
   of $A$ and ${\cal E}$.
 \end{proof}

\begin{remark}\label{remark ideal-involution} Let $A$ be an associative pair and let ${\cal I} $  be a $\ast$-ideal of its  standard imbedding ${\cal E}$.  Then
$ I={\cal I}\cap A=( {\cal
I}_{12},  {\cal
I}_{21})$ is a $\ast$-ideal of $A$, since clearly
 $({\cal I}_{ij})^\ast=(e_i{\cal I}e_j)^\ast= e_j^\ast {\cal I}^\ast e_i^\ast=e_i  {\cal I}e_j= {\cal I}_{ij}$, for  $i\neq j$.
We also note  here  that
$({\cal I}_{11})^\ast={\cal I}_{22}$.
   \end{remark}

\begin{lema} If  $I=(I^+,I^-)$ is  an $\ast$-ideal of an associative pair $A$ with involution $\ast$, and ${\cal I}$ is the ideal of ${\cal E}$ generated by $I$ as in Lemma \ref{lemma ideals 1}(ii), then  ${\cal I}$   is a $\ast$-ideal of the standard imbedding ${\cal E}$
of $A$. Moreover, if $A$ is semiprime and $I$ an essential $\ast$-ideal of $A$, then  ${\cal I}$ is an essential $\ast$-ideal of ${\cal E}$.
\end{lema}

\begin{proof} The first assertion is straightforward, and if $A$ is semiprime and $I$ is essential, the essentiality of ${\cal I}$ follows as in the proof of Lemma \ref{lemma ideals 1}.
\end{proof}

\begin{proposicion} Let $A$ be a semiprime associative pair with involution~$\ast$. Then $(\sum_{i=1}^n \lambda_i a_i^\sigma)^\ast = \sum_{i=1}^n \lambda_i^\ast (a_i^\sigma)^\ast$, where $\lambda_i\in {\cal C}(A)$ and $a_i^\sigma\in A^\sigma$, $\sigma=\pm$, $i=1, \ldots,n$, defines an involution on the central closure  ${\cal C}(A)A$ of $A$, extending the one of $A$.
\end{proposicion}

\begin{proof} Clearly the  involution $\ast$ of $A$, already extended to ${\cal C}(A)$ of $A$  in Proposition \ref{proposition * extended centroid}, also extends to an involution on $ {\cal C}(A) \otimes A$ given by $(\lambda \otimes a^\sigma)^\ast = \lambda^\ast \otimes (a^\sigma)^\ast$, $\sigma=\pm$. Let now $\mu=\overline{(g,I)}\in {\cal C}(A) $ and $y^\sigma\in I^\sigma$,  $\sigma=\pm$. Then we have $(\mu\otimes y^\sigma -1 \otimes g^\sigma( y^\sigma))^\ast = \mu ^\ast\otimes (y^\sigma) ^\ast -1 \otimes (g^\sigma( y^\sigma)) ^\ast  $. Thus since
 $\mu^\ast=\overline{(g^\ast ,I^\ast )} $ with $g^\ast((y^\sigma)^\ast)=( g^\sigma((y^\sigma)^\ast)^\ast)^\ast=( g^\sigma( y^\sigma ))^\ast$, this implies that the ideal $R=(R^+,R^-)$ defined in \ref{subsect central closure AP} is a $\ast$-ideal of $A$. Hence ${\cal C}(A)A$ inherits the involution, also denoted by $\ast$, given by
$(\sum_{i=1}^n \lambda_i a_i^\sigma)^\ast = \sum_{i=1}^n \lambda_i^\ast (a_i^\sigma)^\ast$, where $\lambda_i\in {\cal C}(A)$ and $a_i^\sigma\in A^\sigma$, $\sigma=\pm$, $i=1, \ldots,n$.
\end{proof}

\begin{apartado}\label{subsect ast-central closure AP}  As we did in the previous section, it is possible to define  the scalar extension
$ {\cal C}_\ast(A)A$ of $A$, that will called the
  {\sl $\ast$-central   closure} of $A$. Then $ {\cal C}_\ast(A)A$  is endowed with an involution $(\sum_{i=1}^n \lambda_i a_i^\sigma)^\ast = \sum_{i=1}^n \lambda_i  (a_i^\sigma)^\ast$, where $\lambda_i\in {\cal C}_\ast(A)$ and $a_i^\sigma\in A^\sigma$, $\sigma=\pm$, $i=1, \ldots,n$.
\end{apartado}

\begin{teorema}\label{iso ast central closures} Let $A$ be a semiprime associative pair  with an involution $\ast$ and standard imbedding ${\cal E}$. Then the standard imbedding of $\ast$-central   closure  ${\cal C}_\ast(A)A$ of the associative pair $A$ is isomorphic to
    ${\cal C_\ast(E)E}$.
\end{teorema}

\begin{proof} It follows  from Theorem \ref{th iso  central closures} considering Theorem \ref{th iso ast-centroid}.   \end{proof}

\section{Associative pairs with local PI-algebras }

In this section we extend to associative pairs some of the main results on associative rings with (generalized) polynomial identities, such as the ones due to
Amitsur, Kaplansky, Martindale and Posner Theorems. We refer the reader to \cite{BMM-rings GPI,jaconson PI-algebras,Rowen proceedings,Rowen PI ring}  for quite complete expositions of the classical results of the associative  theory of (generalized) polynomial identities.

\begin{apartado}\label{subsect def strongly primitive AP}
Borrowing the analogous notion  from associative rings  \cite[p.~48]{Rowen proceedings}  (or \cite[p.~281]{Rowen PI ring}), we will say that an associative pair $A$ is {\sl strongly primitive} if $Soc(A)\neq0$, and $A$ is a dense subpair of ${\cal H}=\left(Hom_\Delta(M^-,M^+),Hom_\Delta(M^+,M^-)  \right) $ for a  suitable pair of right vector spaces $M^+$ and $M^-$ over a division PI-ring $\Delta$.
\end{apartado}

\begin{remark}\label{remark AA stronlgy primitive} Strongly primitive associative algebras  are described in \cite[7.5, 7.6]{Rowen PI ring}.  Here limit ourselves to recall that an associative algebra $R$ is strongly primitive if and only if $R$ is primitive and has nonzero PI-ideal. As noted in \cite[1.3]{pi-i} in any strongly primitive  associative algebra the socle and the PI-ideal coincide. (A similar characterization for Jordan systems   was proved  in \cite[Thoerem 4.6]{pi-i}. In \cite{pi-i}  a Jordan system is called rationally primitive (and not strongly primitive) if it is primitive and has a nonzero PI-element. Motivation for that slightly different terminology is given in \cite[4.1]{pi-i}.)
 \end{remark}

\begin{teorema}\label{th strong primitive}
Let $A$  be an associative pair. Then  $A$ is strongly primitive if and only if its standard imbedding ${\cal E}$ is strongly primitive.
\end{teorema}\begin{proof}
Let $A$ be an associative pair, and assume that $A$ is strongly primitive.  Then  $Soc(A)$ is nonzero, and therefore by \ref{subsect socle SI} and Lemma \ref{lemma ideals 1}, the socle $Soc({\cal E})$  of its  standard imbedding ${\cal E}$ is also nonzero. By primitivity of $A$, there are  two  right vector spaces $M^+$ and $M^-$ over a division PI-ring $\Delta$, such that
$A$ is a dense subpair of ${\cal H}=\left(Hom_\Delta(M^-,M^+),Hom_\Delta(M^+,M^-)   \right) $. Then it follows from \cite[2.3]{CGM-primitive AP} that ${\cal E}$ is a primitive associative algebra with   faithful irreducible right ${\cal E}$-module  $M=M^-\oplus M^+$ (over the same division PI-ring $\Delta= End_{{\cal E}}(M^+\oplus M^-)$). Hence ${\cal E}$  is strongly primitive.

Conversely,  let $A$ be an associative pair having a strongly primitive standard imbedding ${\cal E}$.  Now we have  $Soc({\cal E})\neq 0$, hence again   by \ref{subsect socle SI} and Lemma \ref{lemma ideals 1}, $A$ has nonzero socle. Suppose now that ${\cal E}$ is dense in $End(M_\Delta) (=End_\Delta(M))$
 for a right vector space  $M $   over a division PI-ring $\Delta$. As noted in \cite[p. 2598]{CGM-primitive AP}), $M=M^-\oplus M^+$, with $M^+=Me_2$ and $M^-=Me_1$, and then $(M^+,M^-)$ is a faithful irreducible $A$-module over the same division PI-ring $\Delta$. Therefore  $A$ is a strongly primitive associative pair.
\end{proof}

We will apply now this result to obtain an analogue for associative pairs of Amitsur's theorem on primitive algebras with a GPI (see \cite[7.2.9, 7.4.6]{Rowen PI ring}).  Here, as mentioned in the introduction, our GPIs will be nonzero local PI-algebras, so that our version of the GPI condition for a semiprime associative pair $A$ will be the condition $PI(A)\ne 0$. This approach follows the one of \cite{Rowen proceedings}, which in turn relies on the method of ``viewing"  a generalized identity as a polynomial identity of a left (or right) ideal, a method which, according to Rowen \cite[p.~38]{Rowen PI ring}, was initiated by Jain. This was also the approach followed in \cite{pi-ii}, where PI left ideals (which do no make sense in the Jordan theoretical context of that paper) are substituted by PI-elements.

\begin{remark}\label{remark socle AA rowen}  The socle of a primitive ring $R$ can be characterized as the set of all elements of finite rank \cite[Theorem 7.1.13]{Rowen PI ring}. Indeed recall that the  socle $Soc(R)$ of a semiprime ring $R$ is defined  (when  is nonzero)  to be the sum of all  minimal left (equivalently right) ideals of $R$  \cite[Definition~1, Proposition 7.1.6]{Rowen PI ring}, that is the sum of all left (or right) ideals generated by rank one elements \cite[Lemma 7.1.11]{Rowen PI ring}.
 \end{remark}

\begin{apartado}\label{subsect socle rowen} The main ideas given in \cite[p.~254-257]{Rowen PI ring}  can also be applied to
primitive associative pairs leading to
similar results to those mentioned in Remark \ref{remark socle AA rowen} above.  Let $A=(A^+,A^-)$ be a primitive associative pair and suppose that up to isomorphism   $A$ is a dense subpair of a pair $${\cal H}=\left(Hom_\Delta(M^-,M^+),Hom_\Delta(M^+,M^-)   \right),$$ for a faithful  irreducible $A$-module $(M^+,M^-)$.

Then $M=M^-\oplus M^+$ is a faithful irreducible module over the standard imbedding ${\cal E}$ of $A$  \cite[2.3]{CGM-primitive AP}, and $M$  (resp. $(M^+,M^-)$) is a left vector space (resp. a pair of left vector spaces) over the division algebra $\Delta= End(M_{{\cal E}})$. For any associative pair element $a^\sigma\in A^\sigma$, we define the {\sl rank of $a^\sigma$ in $A$} (also the $(M^+,M^-)$-rank of  $a^\sigma$)   to be:
$$rank(a^\sigma)=[M^{-\sigma}a^\sigma:\Delta] = [Ma^\sigma:\Delta], \quad \sigma=\pm.$$
Hence the rank  of $a^\sigma$ is the same independently of whether the element is considered an associative pair element  or an element of the stantard imbedding.

We can add the following remark that links the approach through the standard embedding and the local algebra approach mentioned above, and whose proof is an easy exercise in associative theory: with the notations above, if $a^\sigma \in PI(A^\sigma)$, the local algebra $A^{-\sigma}_{a^{\sigma}}$ is isomorphic to the matrix algebra $M_{t}(\Delta)$ for  $t = rank(a^\sigma)$ (see below the proof of theorem 2.7).
\end{apartado}

\begin{teorema}\label{th AP primitive socle} Let $A$ be a primitive associative pair.  Then, for $\sigma=\pm$:
\begin{enumerate}
\item[(i)] If $a^\sigma\in A^\sigma$ has $rank(a^\sigma)=1$, then $a^\sigma A^{-\sigma} A^\sigma$ is a minimal left ideal of $A$.
\item[(ii)] If $a^\sigma\in A^\sigma$ has $rank(a^\sigma)=t\geq 1$, then there exist rank one elements $a^\sigma_1,\ldots, a^\sigma_t\in A^\sigma$ such that $a^\sigma=\sum_{i=1}^ta^\sigma_i$.
\item[(iii)] $Soc(A^\sigma)$ is the set of elements of finite rank.
\end{enumerate}\end{teorema}

\begin{proof} As mentioned in \ref{subsect socle rowen} above it suffices to review \cite[Lemma~7.1.11, Lemma~7.1.12, Theorem~7.1.13]{Rowen PI ring} introducing the obvious changes to obtain the corresponding results for associative pairs.
\end{proof}

Now we can state the announced analogue of Amitsur's theorem for Associative pairs.

\begin{teorema}\label{th AP GPI  Amitsur}
Let $A$  be an associative pair. Then the following are equivalent:
\begin{enumerate}
\item[(i)] $A$ is strongly primitive.
\item[(ii)] $A$ is  prime and $Soc(A)=PI(A)\neq0$.
\item[(iii)] $A$ is prime and the local algebra at some nonzero element is a simple unital $PI$-algebra.
\end{enumerate}
\end{teorema}\begin{proof}
$(i)\Rightarrow (ii)$ Let $A$ be a strongly primitive associative pair. Then $A$ is prime (see \ref{subsect AP regularity}) and  by Theorem \ref{th strong primitive} its standard imbedding ${\cal E}$ is a strongly primitive algebra. Thus, by \cite[Proposition 7.5.17]{Rowen PI ring}, $Soc({\cal E})=PI({\cal E})$ is a nonzero ideal of ${\cal E}$ and  the equality    $Soc(A)=PI(A)$
follows from \ref{subsect socle SI} and Proposition \ref{prop PI-ideal}. Moreover $Soc(A)=PI(A)\neq 0$  by Lemma \ref{lemma ideals 1}.

$(ii)\Rightarrow (iii)$  Suppose now that $A$ is a prime associative pair with nonzero socle equal to its PI-ideal, and take a nonzero element $0\neq a\in A^{-\sigma}=Soc(A^{-\sigma})=PI(A^{-\sigma})$. By Theorem \ref{th AP primitive socle}(iii)  we can assume $a$ has finite rank $rank(a)=t$. Suppose also that  $a=a^-\in A^-$. Then $A$ is primitive by \cite[2.8]{CGM-primitive AP}, hence a dense subpair of  ${\cal H}=\left(Hom_\Delta(M^-,M^+),Hom_\Delta(M^+,M^-)   \right) $ for a  suitable pair of right vector spaces $M^+$ and $M^-$ over a division ring $\Delta$. Write $M=M^-\oplus M^+$ and take the local algebra $A^+_a$. Then $\widetilde{M} =M/lann_M(a)$ where
$\displaystyle lann_M(a)=\{m\in M=M^-\oplus M^+\mid ma=0\}=M^-\oplus lann_{M^+}(a) $
is a faithful irreducible right $A^+_a$-module, and $[\widetilde{M}: \Delta]=rank (a)=t$. Thus, the local algebra
$A^+_a$ is  a primitive associative PI-algebra satisfying the polynomial standard identity $S_{2t}$. Hence, by Kaplansky's Theorem \cite[Theorem 1.5.16]{Rowen PI ring},    $A^+_a$ is simple unital (since the element $a$ belongs to $Soc(A)$, hence it is  von Neumann regular and the local algebra $A^+_a$  is unital) $PI$-algebra.

$(iii)\Rightarrow (i)$  Consider now a prime associative pair having a simple unital local PI-algebra $A^\sigma_a$  for some $0\neq a\in A^{-\sigma}$. Then, since  as noted in \ref{subsect def local algebra}, the local algebra ${\cal E}_a$ of the standard imbedding ${\cal E}$ of $A$  at the same element satisfies ${\cal E}_a\cong A^\sigma_a$,  we have that ${\cal E}_a$ is a simple unital  associative PI-algebra. Moreover ${\cal E}$  is prime (see \ref{subsect regularity conditions}), hence the standard imbedding ${\cal E}$  of $A$ is strongly primitive by \cite[Proposition 7.5.17(ii)]{Rowen PI ring}. Now the strongly primitivity of the associative pair $A$ follows from
  Theorem \ref{th strong primitive}.
\end{proof}

We next address our Pair analogue of Kaplansky Theorem on PI algebras. Here, as in \cite{pi-ii}, the polynomial identities that we will consider will be homotope polynomial identities, so that we can make use of our results on PI-elements. We begin recalling some facts on homotope polynomial  identities.

 \begin{apartado} Homotope polynomials are the images of associative polynomials $f(x_1,\ldots,x_n)$ of the  free associative algebra $FA[X\cup\{z\}]$    on a countable set of generators $X$ and $z\not\in X$ under  homomorphims $FA[X]\to FA[X\cup\{z\}]^{(z)}$, extending the identity on $X$. Homotope polynomials are usually denoted by  $f(z; x_1,\ldots,x_n)=f(x_1,\ldots,x_n)^{(z)}$.  An associative pair $A$ satisfies a  {\sl homotope PI} (HPI for short), equivalently, $A$ is an {\sl associative HPI-pair } if there exists  $f(x_1,\ldots,x_n)\in FA[X]$
such that $f(y^{-\sigma}; x_1^\sigma,\ldots,x_n^\sigma)$ vanishes under all substitutions of elements $y^{-\sigma}\in A^{-\sigma}$, $ x_i^\sigma
\in A^\sigma$, $\sigma=\pm$. Note that any HPI is in particular a generalized polynomial identity (GPI). Indeed if an associative pair is homotope-PI, then all its local algebras satisfy the same PI.
 \end{apartado}

\begin{teorema}\label{th AP HPI  Kaplansky}
Let $A$  be a primitive associative pair.
\begin{enumerate}
\item[(i)] If the local algebra at any element of $A$ is PI, then $A$ is simple, equal to its socle.
\item[(ii)] If  $A$ is HPI, then $A$ is simple,  equal to its socle.
 \end{enumerate}
\end{teorema}

\begin{proof}  Clearly (ii) is a straightforward  consequence of (i) since all local algebras of
   an associative HPI pair satisfy the same polynomial identity. To prove (i)
  consider $A$ to be a primitive associative pair such that all its local algebras are PI-algebras.  Then, by Theorem \ref{th AP GPI  Amitsur}, we have $A=PI(A)=Soc(A)$, hence $A$ is simple by \cite[Theorem~1]{CGCM-sp AP miin-i}, since it is prime by \ref{subsect AP regularity}, so  \ref{subsect socle AP}  applies.
\end{proof}

\begin{apartado}\label{subsect AP FC} Simple associative pairs coinciding with their socle  are of the form ${\cal F} ({\cal P}_1, {\cal P}_2)$ for two pairs of dual vector spaces
${\cal P}_i=(X_i,Y_i)$, $i=1,2$,
 over the same associative division algebra $\Delta$ \cite[Theorem 2]{CGCM-sp AP miin-i}.   Recall, see \cite[Remark p.~483]{CGCM-sp AP miin-i}, that
  such an  associative pair has finite capacity if and only if one of the
 vector spaces is finite-dimensional over $\Delta$.
 The classification of simple associative pairs having finite capacity is given in \cite[Theorem 11.16]{loos-jp}.
\end{apartado}

The existence of a homotope-PI on a primitive
 associative pair provides a bound for the dimension of (at least) one of the pairs of dual vector spaces, ensuring then the finite capacity of the associative pair.

\begin{corolario}\label{corollary - th AP HPI Kaplansly}
Let $A$  be a primitive associative pair. If $A$ satisfies a  homotope-PI of degree $d$, then $A$ has finite capacity at most $ [d/2]$.
\end{corolario}

\begin{proof}   Assume that $A$ is a primitive associative pair satisfying a  homotope-PI of degree~$d$. Then, by Theorem \ref{th AP HPI  Kaplansky}, $A$  is simple and equal to its socle. Indeed, by \cite[Theorem 2]{CGM-primitive AP},  $A=Soc(A)=({\cal F}_\Delta(M^-,M^+)  ,{\cal F}_\Delta(M^+,M^-) )$, where $(M^+,M^-)$ is a faithful irreducible right $A$-module,  and $\Delta=End_{{\cal E}}(M^+\oplus M^-)$  is a division algebra where ${\cal E}$ is  the standard imbedding of $A$.

Take now an element $a\in A^{-\sigma}=Soc(A^{-\sigma})$. By Theorem \ref{th AP primitive socle}(iii)  we can assume that $a$ has  finite rank $ rank (a)=r$. Then the local algebra
$A^\sigma_a$ is simple (\ref{subsect def local algebra}), and according to
  \cite{pares-PI} (see also \cite{paniello-tesis})  we can assume  that
$A^\sigma_a$ is contained into a matrix algebra $M_r(F)$, where $F$ denotes a maximal subfield of the division ring $\Delta$.
As a result, by \cite[1.4.1]{Rowen PI ring} $A^\sigma_a$  satisfies the standard identity $S_{2r}$, hence $2r\leq d$. Thus $A$ has finite capacity at most  $ [d/2]$.
\end{proof}

Next, the  closeness between the central  closures of semiprime associative pairs and that of their standard  imbeddings makes it possible to obtain the following associative pair version of Martindale Theorem for prime associative algebras satisfying a generalized polynomial identity.

\begin{teorema}\label{th AP GPI  Martindale}   Let $A$ be a prime associative pair. If $PI(A)\neq 0$, then the central  closure
${\cal C}(A)A$ of $A$ is a primitive associative pair with nonzero socle equal to $PI({\cal C}(A)A)$. Moreover $PI(A)=A\cap Soc({\cal C}(A)A)$.
\end{teorema}

\begin{proof} Let $A$ be a prime associative pair having nonzero  PI-elements. Then, by \ref{subsect regularity conditions},
its standard imbedding ${\cal E}$ is a prime   associative algebra with nonzero PI-ideal
$PI({\cal E})$ by Proposition \ref{prop PI-ideal}.  Hence the central closure  ${\cal C(E)E}$  of ${\cal E}$   is a strongly primitive associative algebra with nonzero socle   by
\cite[Theorem 7.6.15]{Rowen PI ring}. Moreover, since  by \cite[Proposition 7.5.17]{Rowen PI ring},   $Soc({\cal C(E)E})=PI({\cal C(E)E})$, it holds that   $PI({\cal E})=  {\cal E}\cap  Soc({\cal C(E)E})$ .

On the other hand, by Theorem
\ref{th iso central closures},  we
  have that ${\cal C(E)E}$ is isomorphic to the standard imbedding of the central closure ${\cal C}(A)A$ of the associative pair $A$. Thus, by Theorem \ref{th strong primitive}, ${\cal C}(A)A$ is a strongly primitive associative pair with $Soc({\cal C}(A)A)= PI({\cal C}(A)A)\neq 0$ as a result of Theorem \ref{th AP GPI  Amitsur}.

  Finally we claim that
     $PI(A)=A \cap PI({\cal C}(A)A)=A\cap Soc({\cal C}(A)A)$. Note that it suffices to prove  $PI(A) \subseteq A \cap PI({\cal C}(A)A) $.
     Take $a=a^{-\sigma}\in PI(A^{-\sigma})$. Then, by Proposition \ref{prop PI-ideal}, $a\in A^{-\sigma} \cap PI( {\cal E}) =A^{-\sigma}  \cap Soc({\cal C(E)E})$, which implies that  $a$ is von Neumann regular in $ {\cal C(E)E} $.  Then, as a result of
     \ref{subsect def local algebra}, Theorem  \ref{th iso  central closures}   and Remark \ref{remark Martindale AQ}, we obtain the following contaiments of local algebras  $A^\sigma_a\subseteq {\cal C}(A)A^\sigma_a\cong {\cal C(E)E}_a \subseteq Q_s({\cal E} )_a$, where $Q_s({\cal E} )$ denotes the Martindale symmetric ring of quotients of the standard imbedding ${\cal E}$ of $A$.
   Moreover,  from the   (von Neumann) regularity of $a$ in $ {\cal C(E)E} $, hence in $Q_s({\cal E} )$, we have $Q_s({\cal E})_a\cong Q_s({\cal E}_a)$ by \cite[Theorem 3]{gs-Local AQ}.  Therefore $A^\sigma_a\subseteq {\cal C}(A)A^\sigma_a\cong {\cal C(E)E}_a \subseteq Q_s({\cal E} )_a\cong Q_s({\cal E}_a)$. Thus, as $A_a^\sigma$ is prime by \ref{subsect def local algebra},  $Q_s({\cal E}_a)$ is a PI-algebra by 
   \cite[Corollary 6.1.7]{BMM-rings GPI}, and that implies that ${\cal C}(A)A^\sigma_a$ is a PI-algebra, and therefore we finally obtain $a\in A \cap PI({\cal C}(A)A)$.
\end{proof}

\begin{teorema}\label{th AP GPI Posner}
Let $A$  be a  prime associative pair.
\begin{enumerate}
\item[(i)] If the local algebra    $A_a^{-\sigma}$ at each element $a\in A^{-\sigma}$ of $A$ is PI, then the  central closure  ${\cal C}(A)A$ of $A$ is simple equal to its socle.
\item[(ii)] If  $A$ is HPI, then ${\cal C}(A)A$ is simple,  equal to its socle.
 \end{enumerate}
\end{teorema}

\begin{proof}
Under any of the above assumptions, the central  closure  ${\cal C}(A)A$ of $A$ is a strongly primitive associative pair  by Theorem \ref{th AP GPI  Martindale}. Besides, $A=PI(A)$ is contained in $Soc({\cal C}(A)A ) $. Thus   $ {\cal C}(A)A   = Soc({\cal C}(A)A ) $ and the simplicity of ${\cal C}(A)A$  follows from \cite[Theorem 1]{CGCM-sp AP miin-i}.
\end{proof}

As a consequence of the previous results, we obtain the following Associative Pair version of Posner's Theorem:

\begin{corolario}\label{corrollary - th AP GPI Posner}
Let $A$  be a  prime associative pair. If $A$ satisfies a homotope-PI of degree $d$, then its central  closure  ${\cal C}(A)A$  has finite capacity at most  $ [d/2]$.
\end{corolario}

\begin{proof} It follows from Theorem \ref{th AP GPI Posner} as a result of Corollary \ref{corollary - th AP HPI Kaplansly}.
\end{proof}

%
%
 \addcontentsline{toc}{chapter}{\protect\numberline{}  Bibliograf\'{\i}a.}


\end{document}